\documentclass{siamltex}

\usepackage{amsmath, amsfonts, amssymb, latexsym}
\usepackage{graphicx}

\title{A multi-mesh finite element method for Lagrange elements of
  arbitrary degree}

\author{{\sc Axel Voigt}\thanks{Institut f\"ur Wissenschaftliches
    Rechnen, TU Dresden, Zellescher Weg 12-14, 01062 Dresden, Germany
    ({\tt axel.voigt@tu-dresden.de})} \and {\sc Thomas
    Witkowski}\thanks{Institut f\"ur Wissenschaftliches Rechnen, TU
    Dresden, Zellescher Weg 12-14, 01062 Dresden, Germany ({\tt
      thomas.witkowski@tu-dresden.de})}}

\begin{document}
\maketitle

\begin{abstract}
  We consider within a finite element approach the usage of different
  adaptively refined meshes for different variables in systems of
  nonlinear, time-depended PDEs. To resolve different solution
  behaviours of these variables, the meshes can be independently
  adapted. The resulting linear systems are usually much smaller, when
  compared to the usage of a single mesh, and the overall
  computational runtime can be more than halved in such cases. Our
  multi-mesh method works for Lagrange finite elements of arbitrary
  degree and is independent of the spatial dimension.  The approach is
  well defined, and can be implemented in existing adaptive finite
  element codes with minimal effort. We show computational examples in
  2D and 3D ranging from dendritic growth to solid-solid
  phase-transitions.  A further application comes from fluid dynamics
  where we demonstrate the applicability of the approach for solving
  the incompressible Navier-Stokes equations with Lagrange finite
  elements of the same order for velocity and pressure. The approach
  thus provides an easy to implement alternative to stabilized finite
  element schemes, if Lagrange finite elements of the same order are
  required.
\end{abstract}

\section{Introduction}
Nowadays, adaptive {\em h}-methods for mesh refinement are a standard
technique in finite element codes. They are used to resolve a mesh due
to the local behaviour of the solution. When solving PDEs with
multiple variables, e.g.\ velocity and pressure in the Navier-Stokes
equations, the mesh has to be adapted to the behaviour of all
components of the solution. If these behaviours are different, the use
of a single mesh may lead to an inefficient numerical method. In this
work we propose a multi-mesh finite element method that makes it
possible to resolve the local nature of different components
independently of each other. This method works for Lagrange elements
of arbitrary degree in any dimension. Furthermore, the method works
``on the top'' of standard adaptive finite element methods. Hence,
only small changes are required if our method has to be implemented in
existing finite element codes. We have implemented the multi-mesh
method in our finite element software AMDiS (adaptive multidimensional
simulations), see \cite{Vey_2007}, for Lagrange finite elements up to
fourth degree for 1D, 2D and 3D.

To our best knowledge, Schmidt \cite{Schmidt_2003} was the first who
has considered the use of multiple meshes in this context. Li et al.\
have introduced a very similar technique and used it to simulate
dendritic growth \cite{Li_2005, DiLi_2009, HuLiTang_2009}. Although
introducing a multi-mesh technique, in none of these publications the
method is formally derived.  Furthermore implementation issues are not
discussed and detailed runtime results, which compare the overall
runtime between the single-mesh and the multi-mesh method are
missing. In contrast, in this work we will formally show how multiple
meshes are used in the context of assembling matrices and vectors in
the assembly step of the finite element methods and will discuss
issues related to error estimates for each component. Furthermore, we
will compare the runtimes of both methods and show that the multi-mesh
method is superior to the single-mesh method, when one component of
the PDE can locally be resolved on a coarser mesh. Solin et al.\
\cite{Solin_2007, Solin_2010} have also introduced a multi-mesh
method, but for {\em hp}-FEM. Their method is based on transforming
quadrature points which is harder to implement in existing finite
element codes. Furthermore, also in these works detailed runtime
studies are missing.  We should further mention other approaches which
are commonly used to deal with different meshes for different
components of coupled systems. Especially in the case of multi-physics
applications a need exists to couple independent simulations code. A
standard tool which can be used to couple various finite element codes
is MpCCI (mesh-based parallel code coupling interface)
\cite{Joppichetal_CCPE_2006}.  In this approach an interpolation
between the different solutions from one mesh to the other is
performed which for different resolutions of the involved meshes will
lead to a loss in information and is thus not the method of choice for
the problems to be discussed in this work.

The paper is structured as follows: In the next section we give a
brief overview on adaptive meshes, especially in the context of AMDiS,
and introduce the terminology used throughout this paper. Section 3
introduces the so called {\em virtual mesh assembling}, which is the
basis of our multi-mesh method. It is shown, how the coupling meshes
are build in a virtual way and how the corresponding coupling
operators are assembled on them. In Section 4, we present several
numerical experiments in 2D and 3D that show the advantages of the
multi-mesh method. The last section summarizes our results.

\section{Adaptive meshes}
The usage of an adaptive mesh, together with error estimators and a
refinement strategy, is a standard finite element technique to compute
solution of PDEs with a given accuracy with the lowest possible
computational effort. For a general overview on this topic see for
example \cite{Verfuert_Teubner_1996}, and references
therein. In this section, we describe how adaptive meshes and the
associated algorithms are implemented in our finite element software
AMDiS. The data structures that are used to store and manipulate
adaptive meshes, are the basis for a fast and efficient multi-mesh
method, as it is presented in the next section.

\begin{figure}
  \centering
  \includegraphics[scale=0.8]{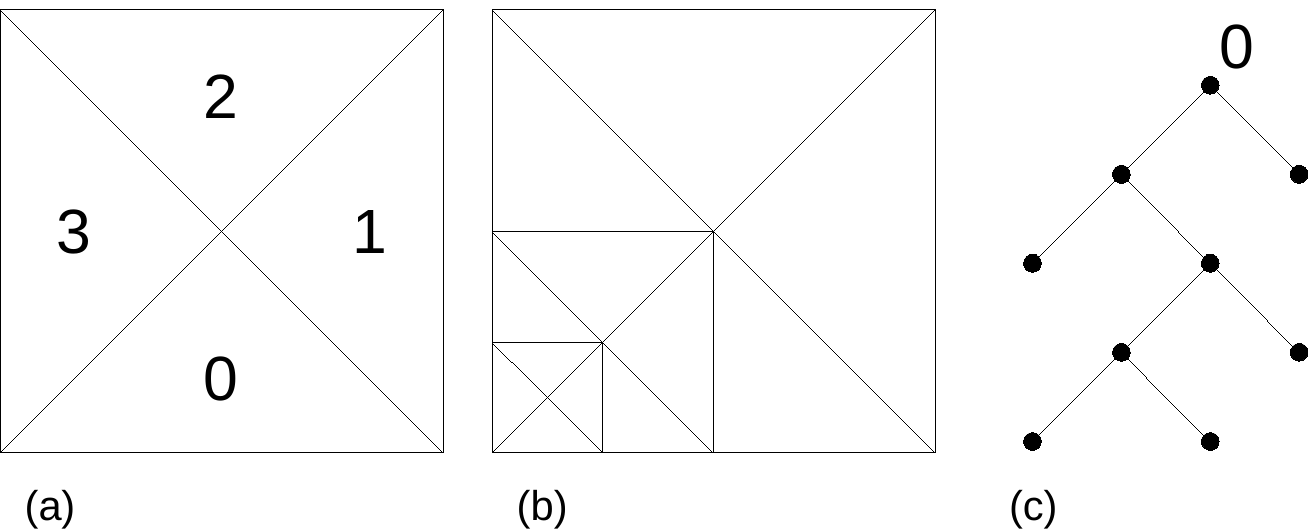}
  \caption{(a) A two-dimensional macro mesh. (b) Some refinements of
    it. (c) Binary tree of macro element 0}
  \label{img:adaptmesh}
\end{figure}

A mesh in AMDiS consist of simplicial elements, which are lines in 1D,
triangles in 2D and tetrahedral in 3D. If an element has to be refined
within the adaption loop, it will be bisected into two elements of the
same dimension. The refinement algorithm, that is implemented in
AMDiS, is described in \cite{Schmidt_Siebert_2005} in more detail. The
two new elements are called {\em children} of the {\em parent
  element}. The coarsest elements are called {\em macro
  elements}. Accordingly, the {\em macro mesh} is the union of all
macro elements. The refinement history of a macro element is
represented by a binary tree. If a node in this tree is a leaf, the
corresponding element is not refined. Otherwise, the children of the
node represent the children of the element. To avoid hanging nodes
(vertices of an element which are not vertices of the neighbour
element), it may be necessary to first refine a neighbour of an
element before refining the element itself. In this way, a single
refinement can cause some propagation refinements of elements in the
neighbourhood. In Figure \ref{img:adaptmesh}, a triangular macro mesh
consisting of four elements (a), some refinements of this macro
triangulation (b) and the corresponding binary tree for macro element
0 (c) are shown.

All topological and geometrical information are only stored for the
macro elements. To get these information for other elements in the
mesh, an algorithm called {\em mesh traverse} is used. This algorithm 
traverses all binary trees, corresponding to the macro
elements, and recursively computes the requested information, e.g.\
coordinates of vertices, for the requested child elements from the
information of its parent element. The mesh traverse algorithm than
calls a function that computes on the element data. An algorithm, that
makes use of the mesh traverse, can specify the element level on which
it wants to process. The {\em level} of an element is defined to be
the depth of its node representation in the binary tree. In most
cases, the {\em leaf level}, i.e., the union of all leaves in all
binary trees, is used for mesh traverse.

The additional effort to compute the element data for the leaf mesh
multiple times is usually less than 1\% of the computational time of
the function that process on the mesh. Instead, a large amount of
memory can be saved. In AMDiS, the information data for one element
requires around 200 bytes per element. Many of our simulations make
use of around one million elements per processor. Storing all element
information explicitly for these simulations would require more than
200 Mbyte of extra memory.

\subsection{Error estimation and adaptive strategies}

AMDiS provides different methods for error estimation and mesh
refinement. Furthermore, due to an abstract interface, is is easily
possible to implement other methods, that are more adapted to a
specific PDE. In the following, we describe the most common ones.

In AMDiS, the standard estimator for the spatial error is the residual
error estimator as for example described by Verf\"urth
\cite{Verfurth_1994a} for general non-linear elliptic PDEs. This local
element wise estimator defines for each element $T$ of a mesh an
indicator, depending on the finite element solution $u_{h}$, by
\begin{equation}
  \label{eq:residual_ee}
  \eta_{T}(u_{h}) = C_{o} R_{T}(u_{h}) + C_{1} \sum_{E \subset T} J_{E}(u_{h}),
\end{equation}
where $R_{T}$ is the element residual on element $T$ and $J_{E}$ is
the jump residual, defined on an edges $E$. $C_{0}$ and $C_{1}$ are
user defined constants that can be used to adjust the estimator. The
global error estimation of a finite element solution is than the sum
of the local estimates:
\begin{equation}
  \eta(u_{h}) = \left( \sum_{T \in {\cal T}} \eta_{T}(u_{h})^{p} \right)^{1/p}, \hspace*{1.0cm} p \geq 1,
\end{equation}
where ${\cal T}$ is a partitions of the domain into simplices. The
user has to specify an error tolerance $tol$. If $\eta(u_{h}) > tol$,
the mesh must be refined in some way in order to reduce the
error. Several strategies are implemented in AMDiS:
\begin{itemize}
\item the maximum strategy, as described by Verf\"urth
  \cite{Verfurth_1994b},
\item the equidistribution strategy, as described by Eriksson and
  Johnson \cite{Eriksson_1991},
\item the guaranteed error reduction strategy, as described by
  D\"orfler \cite{Dorfler_1996}.
\end{itemize}
In all of the numerical experiments in Section
\ref{sec:numerical_results}, we make use of the equidistribution
strategy, which we found out to be easy to use and which provided good
results. The basis for this strategy is the idea that the error is
assumed to be equidistributed over all elements. Hence, if the
partition consists of $n$ elements, the element error estimates should
fulfill
\begin{equation}
  \eta_{T}(u_{h}) \approx \frac{tol}{n^{1/p}} \equiv \eta_{eq}(u_{h}).
\end{equation}
The equidistribution strategy now introduces two parameters,
$\theta_{R} \in (0, 1)$ and $\theta_{C} \in (0, 1)$. An element is
refined, if $\eta_{T}(u_{h}) > \theta_{R} \eta_{eq}(u_{h})$, and the
element is coarsen, if $\eta_{T}(u_{h}) \leq \theta_{C}
\eta_{eq}(u_{h})$. The parameter $\theta_{R}$ is usually chosen to be
close to 1, and $\theta_{C}$ to be close to 0, respectively.

If the equation consists of more than one variable, one can define
multiple error estimators and mesh adaption strategies. This is also
the case, if only one mesh is used to resolve all variables. The error
estimators work independently of each other, and the user can provide
the constants $C_{0}$ and $C_{1}$ for each component. Correspondingly,
independent mesh adaption strategies may be defined for each
component. In the case of the single-mesh method, an element is
refined, if at least one strategy has marked the element to be
refined. An element is coarsen, if all strategies have marked it to be
coarsen. In the multi-mesh method, these conditions are omitted,
because the meshes are adapted independently of each other.

\section{Virtual mesh assembling}
The basis of our multi-mesh method is the so called {\em virtual mesh
  assembling}. Systems of PDEs usually involve coupling terms. If each
component of the system is assembled on a different mesh, special care
has to be devoted to these coupling terms. In the next section, we
shortly describe this situation. The section
\ref{sec:dual_mesh_traverse} than introduces the {\em dual mesh
  traverse}.  This algorithm create a virtual union of two meshes
without creating it explicitly. To the last, we show how to compute
integrals, that appear within the assemble procedure, on these virtual
meshes.

\subsection{Coupling terms in systems of PDEs}

To illustrate the techniques presented in this section, we consider
the homogeneous biharmonic equation as a simple example for general
systems of PDEs. This equation reads:
\begin{equation}
  \Delta^2 u = 0 \;\; \text{in} \;\; \Omega \hspace*{0.5cm} \text{and} \hspace*{0.5cm} 
  u = \frac{\partial u}{\partial n} = 0 \;\; \text{on} \;\; \partial \Omega,
\end{equation}
with $u \in C^{4}(\Omega) \: \cap \: C^{1}(\bar{\Omega})$. Using
operator splitting, the biharmonic equation can be rewritten as a
system of two second order elliptic PDEs:
\begin{equation}
\begin{split}
-\Delta u + v &= 0 \\
\Delta v &= 0
\end{split}
\end{equation}
The standard mixed variational formulation of this system is: Find
$(u, v) \in H_{0}^{1}(\Omega) \times H^{1}(\Omega)$ such that
\begin{equation}
  \label{eq:biharmonic_var1}
  \begin{split}
    \int_{\Omega}\nabla u \nabla \phi dx + \int_{\Omega} v \phi dx &= 0 \hspace*{1cm} \forall \phi \in H^{1}(\Omega) \\
    \int_{\Omega}\nabla v \nabla \psi dx &= 0 \hspace*{1cm} \forall \psi \in H_{0}^{1}(\Omega)
  \end{split}
\end{equation}
To discretize these equations, we assume that ${\cal T}_{h}^{0}$ and
${\cal T}_{h}^{1}$ are different partitions of the domain $\Omega$
into simplices. Than, $V_{h}^{0} = \{v_{h} \in H^{1} : v_{h}|_{T} \in
P^{n} \;\forall T \in {\cal T}_{h}^{0}\}$ and $V_{h}^{1} = \{v_{h} \in
H_{0}^{1} : v_{h}|_{T} \in P^{n} \;\forall T \in {\cal T}_{h}^{1}\}$
are finite element spaces of globally continuous, piecewise polynomial
functions of an arbitrary but fixed degree. We thus obtain: Find
$(u_{h}, v_{h}) \in V_{h}^{0} \times V_{h}^{1}$ such that
\begin{equation}
  \label{eq:biharmonic_var2}
  \begin{split}
    \int_{\Omega}\nabla u_{h} \nabla \phi dx + \int_{\Omega} v_{h} \phi dx &= 0 \hspace*{1cm} \forall \phi \in V_{h}^{0}(\Omega) \\
    \int_{\Omega}\nabla v_{h} \nabla \psi dx &= 0 \hspace*{1cm} \forall \psi \in V_{h}^{1}(\Omega).
  \end{split}
\end{equation}

Let define $\{\phi_{i} \mid 1 \leq i \leq n\}$ and $\{\psi_{i} \mid 1
\leq i \leq m\}$ to be the nodal basis of $V_{h}^{0}$ and $V_{h}^{1}$,
respectively. Hence, $u_{h}$ and $v_{h}$ can be written by the linear
combinations $u_{h} = \sum_{i = 1}^{n} u_{i} \phi_{i}$ and $v_{h} =
\sum_{i = 1}^{m} v_{i} \psi_{i}$, with $u_{i}$ and $v_{i}$ the unknown
real coefficients. Using these relations and braking up the domain in
the partitions of $\Omega$, Equation (\ref{eq:biharmonic_var2})
rewrites to
\begin{equation}
  \label{eq:biharmonic_var3}
  \begin{split}
    \sum_{j = 1}^{n} u_{j} 
    \left( \sum_{T \in {\cal T}^{0}_{h}} \int_{T} \nabla \phi_{j} \cdot \nabla \phi_{i} \right)
    + \sum_{j = 1}^{m} v_{j} 
    \left( \sum_{T \in {\cal T}^{0}_{h} \cup {\cal T}^{1}_{h}} \int_{T} \psi_{j} \phi_{i} \right)
    &= 0
    \hspace*{0.75cm} i = 1 \dots n \\
    \sum_{j = 1}^{m} v_{j} 
    \left( \sum_{T \in {\cal T}^{1}_{h}} \int_{T} \nabla \psi_{j} \cdot \nabla \psi_{i} \right)
      &= 0
    \hspace*{0.75cm} i = 1 \dots m. \\
  \end{split}
\end{equation}
To compute the coupling term, we have to define the union of two
different partitions $T_{h}^{0} \cup T_{h}^{1}$. For this, we make a
restriction on the partitions: Any element $T^{0} \in {\cal T}_h^{0}$
is either a subelement of an element $T^{1} \in {\cal T}_h^{1}$, or
vice versa. This restriction is not very strict. It is always
fulfilled for the standard bisectioning algorithm, if the initial
meshes for all components share the same macro mesh. Then, ${\cal
  T}^{0}_{h} \cup {\cal T}^{1}_{h}$ is the union of the locally finest
simplices.

The most common way to compute the integrals in Equation
(\ref{eq:biharmonic_var3}) is to define local basis functions. We
define $\phi_{i,j}$ to be the $j$-th local basis function on an
element $T_{i} \in {\cal T}^{0}_{h}$. $\psi_{i,j}$ is defined in the
same way for elements in the partition ${\cal T}^{1}_{h}$.

Because the global basis functions $\phi_{i}$ and $\psi_{j}$ are
defined on different triangulations of the same domain, it is not
straightforward to calculate the coupling term $\int_{\Omega} \psi_{j}
\phi_{i}$ in an efficient manner. For evaluating this integral, two
different cases may occur: either the integral has to be evaluated on
an element from the partition ${\cal T}^{0}_{h}$ or on an element from
${\cal T}^{1}_{h}$. For what follows, we fix the first case. In terms
of local basis functions we have to evaluate
\begin{equation}
  \label{eq:eval_coupling}
  \int_{T_{i} \in {\cal T}^{0}_{h}} \psi_{k,l} \phi_{i,j}
\end{equation}
for some $j$ and $l$, and there exists an element $T_{k} \in {\cal
  T}^{1}_{h}$, with $T_{i} \subset T_{k}$. Our aim is to develop a
multi-mesh method that works on the top of existing finite element
software. These have usually implemented special methods to evaluate
local basis functions, or the multiplication of two basis functions,
respectively, on elements very fast. This involves for example
precalculated integral tables or fast quadrature rules. All these
methods cannot directly be applied to the coupling terms, because
$\psi_{k,l}$ in Equation (\ref{eq:eval_coupling}) is not a local basis
function of the element $T_{i}$. The general idea to overcome this
problem is to define the basis functions $\psi_{k,l}$ by a linear
combination of local basis functions of $T_{i}$. Thus,
\begin{equation}
  \label{eq:linearcombination1}
  \int_{T_{i} \in {\cal T}^{0}_{h}} \psi_{k,l} \phi_{i,j} = 
  \int_{T_{i} \in {\cal T}^{0}_{h}} \sum_{m} (c_{k,m} \phi_{i,m}) \phi_{i,j},
\end{equation}
with some real coefficients $c_{k,m}$. For the other case, i.e., the
integral in the coupling term is evaluated on an element $T_{i} \in
{\cal T}^{1}_{h}$, we have
\begin{equation}
  \label{eq:linearcombination2}
  \int_{T_{i} \in {\cal T}^{1}_{h}} \psi_{k,l} \phi_{i,j} = 
  \int_{T_{i} \in {\cal T}^{1}_{h}} \psi_{k,l} \sum_{m} (c_{i,m} \psi_{k,m}).
\end{equation}

Summarized, to make evaluate of the coupling terms possible, two
different techniques have to be defined and implemented: Firstly, the
method requires to build a union of two meshes. This leads to an
algorithm which we name {\em dual mesh traverse}. It will be discussed
in the next section. Once the union is obtained, we need to calculate
the coefficients $c_{i,j}$ and to incorporate them in the finite
element assemblage procedure such that the overall change of the
standard method is as small as possible.

\subsection{Creation of the virtual mesh}
\label{sec:dual_mesh_traverse}

The simplest way to obtain the union of two meshes is to employ the
data structure they are stored in. Hence, in our case we could
explicitly build the union by joining the binary trees of both meshes
into a set of new binary trees. Especially when we consider meshes
that change in time, this procedure is not only too time-consuming but
requires also additional memory to store the joined mesh. To avoid
this, we exploit the case that in AMDiS functions, that process
element data, never directly work on the mesh data but instead use the
mesh traverse algorithm, that creates the requested element data on
demand. According to this method, we define the {\em dual mesh
  traverse} that traverses two meshes in parallel and thus creates the
union of both meshes in a virtual way.

\begin{figure}
  \centering
  \includegraphics[scale=0.65]{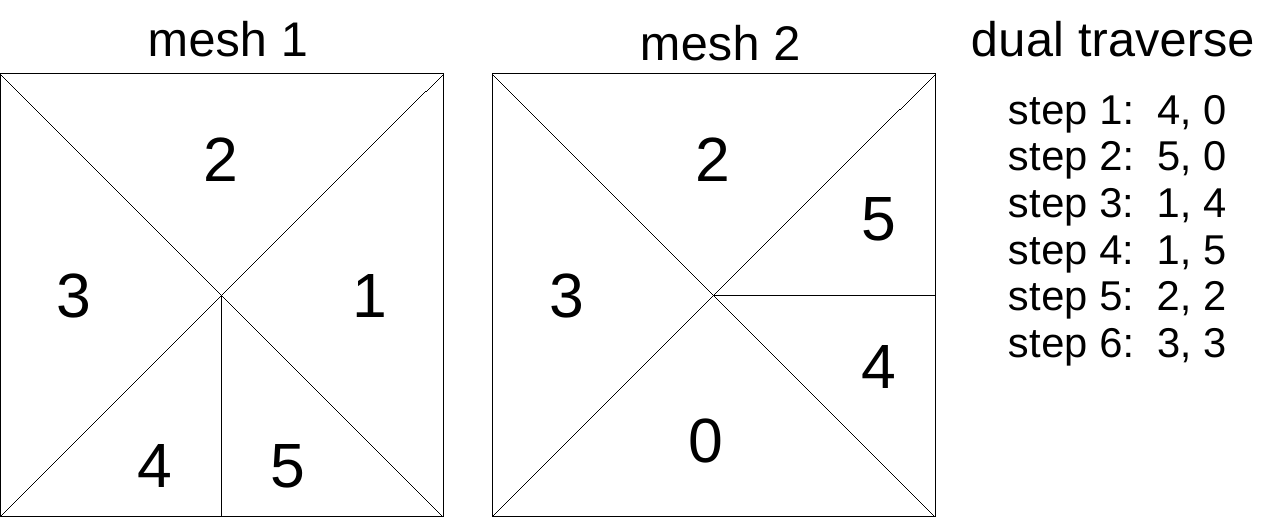}
  \caption{Dual mesh traverse on two independent refined meshes on the
    same macro mesh.}
  \label{fig:dualmesh_traverse}
\end{figure}

For the dual mesh traverse the only requirement is that both meshes
must share the same macro mesh, but they can be refined independently
of each other. Due to this requirement and because of the bisectioning
refinement algorithm the following holds: If the intersection of two
elements of two different meshes is non empty, then either both
elements are equal or one element is a real subelement of the
other. To receive the leaf level of the virtual mesh, the dual mesh
traverse simultaneously traverses two binary trees, each corresponding
the the same macro element in both macro meshes. The algorithm than
calls a user defined function, e.g. the element assembling function or
an element wise error estimator, that works on pairs of elements, with
both, the larger and the smaller element of the current traverse. The
larger of both elements is fixed as long as all smaller subelements in
the other mesh are traversed. Figure \ref{fig:dualmesh_traverse} shows
a simple example for a macro mesh consisting of four macro
elements. In the first mesh macro element 0, and in the second mesh
macro element 1 are refined once.

\subsection{Assembling of element matrices}

To make the overall calculation as fast as possible, the form as given
by the Equations (\ref{eq:linearcombination1}) and
(\ref{eq:linearcombination2}) is not appropriate. To implement these
transformations, some changes of the inner assemblage procedure would
be required. In the following, we show in which way the transform
matrices can be incorporated in the standard assemble procedure, such
that the changes will be as small as possible. To the first, we have
to distinguish two cases: the smaller of both elements defines either
the space of test functions, or it defines the space of trial
functions. For the first case, assume a general zero order term of the
form $\int_{\Omega} \psi_{i} c \phi_{j}$, with $c \in
L^{\infty}(\Omega)$ and some local basis functions $\psi_{i}$ and
$\phi_{j}$, has to be assembled on a virtual mesh. Then, for some
elements $T$ and $T'$, with $T' \subset T$, the element matrix
$M_{T'}$ is given by
\begin{equation}
  \label{eq:zot_matrix_trans}
  \begin{split}
    M_{T'} &=
    \begin{pmatrix}
      \int_{T'} \psi_{0} c \phi_{0} & \dots & \int_{T'} \psi_{0} c \phi_{n} \\
      \vdots & & \vdots \\
      \int_{T'} \psi_{n} c \phi_{0} & \dots & \int_{T'} \psi_{n} c \phi_{n}
    \end{pmatrix} \\
    &= 
    \begin{pmatrix}
      \int_{T'} \sum_{i} (c_{0i}\phi_{i}) c \phi_{0} & \dots & \int_{T'}  \sum_{i} (c_{0i}\phi_{i}) c \phi_{n} \\
      \vdots & & \vdots \\
      \int_{T'} \sum_{i} (c_{ni}\phi_{i}) c \phi_{0}& \dots & \int_{T'} \sum_{i} (c_{ni}\phi_{i}) c \phi_{n}
    \end{pmatrix} \\
    &=
    \begin{pmatrix}
      \sum_{i}c_{0i} \int_{T'} \phi_{i} c \phi_{0} & \dots & \sum_{i}c_{0i} \int_{T'} \phi_{i} c \phi_{n} \\
      \vdots & & \vdots \\
      \sum_{i}c_{ni} \int_{T'} \phi_{i} c \phi_{0} & \dots & \sum_{i}c_{ni} \int_{T'} \phi_{i} c \phi_{n}
    \end{pmatrix} \\
    &=
    {\cal C} \cdot
    \begin{pmatrix}
      \int_{T'} \phi_{0} c \phi_{0} & \dots & \int_{T'} \phi_{0} c \phi_{n} \\
      \vdots & & \vdots \\
      \int_{T'} \phi_{n} c \phi_{0} & \dots & \int_{T'} \phi_{n} c \phi_{n}
    \end{pmatrix},
  \end{split}
\end{equation}
where $\phi_{i}$ are the local basis function defined on $T'$,
$\psi_{i}$ are the local basis function defined on $T$, and ${\cal C}$
is the transformation matrix for the local basis function from $T$ to
$T'$. This shows, that to assemble the element matrix of a virtual
element, there is no need for larges changes within the assemble
procedure. The finite element code needs only to assemble the element
matrix of the smaller element $T'$ and multiply the result with the
transformation matrix. Hence, if the transformation matrices can be
computed easily, the overhead for virtual element assembling can be
neglected.

The second case, where the smaller element defines the space of trial
functions, is very similar. The same calculation as above shows that
the following holds:
\begin{equation}
  \label{eq:zot_matrix_trans2}
  M_{T'} =
  \begin{pmatrix}
    \int_{T'} \phi_{0} c \phi_{0} & \dots & \int_{{\cal T}'} \phi_{0} c \phi_{n} \\
    \vdots & & \vdots \\
    \int_{T'} \phi_{n} c \phi_{0} & \dots & \int_{T'} \phi_{n} c \phi_{n}
  \end{pmatrix}
  \cdot {\cal C}^{T}
\end{equation}

In a similar way we can also reformulate the element matrices for
general first and second order terms. For a general second order term
of the form $\int_{\Omega} \nabla \psi_{i} \cdot {\bf A} \nabla
\phi_{j}$, with ${\bf A}: \Omega \mapsto \mathbb{R}^{d \times d}$, the
element matrix $M_{T'}$ can be rewritten in the same way as we have
done it in (\ref{eq:zot_matrix_trans}) for general zero order terms: 
\begin{equation}
  \begin{split}
    M_{{\cal T}'} &=
    \begin{pmatrix}
      \int_{{\cal T}'} \nabla \psi_{0} \cdot {\bf A} \nabla \phi_{0} & \dots & 
      \int_{{\cal T}'} \nabla \psi_{0} \cdot {\bf A} \nabla \phi_{n} \\
      \vdots & & \vdots \\
      \int_{{\cal T}'} \nabla \psi_{n} \cdot {\bf A} \nabla \phi_{0} & \dots & 
      \int_{{\cal T}'} \nabla \psi_{n} \cdot {\bf A} \nabla \phi_{n}
    \end{pmatrix} \\
    &=
    {\cal C}_{\nabla} \cdot 
    \begin{pmatrix}
      \int_{{\cal T}'} \nabla \phi_{0} \cdot {\bf A} \nabla \phi_{0} & \dots & 
      \int_{{\cal T}'} \nabla \phi_{0} \cdot {\bf A} \nabla \phi_{n} \\
      \vdots & & \vdots \\
      \int_{{\cal T}'} \nabla \phi_{n} \cdot {\bf A} \nabla \phi_{0} & \dots & 
      \int_{{\cal T}'} \nabla \phi_{n} \cdot {\bf A} \nabla \phi_{n}
    \end{pmatrix} \\
  \end{split}
\end{equation}
where the coefficients of the matrix ${\cal C}_{\nabla}$ are defined such that
\begin{equation}
  \label{eq:coef_grad}
  \nabla \psi_{i} = \sum_{j=0}^{n} c_{ij} \nabla \phi_{i}.
\end{equation}
If the smaller element defines the space of trial functions, we can
establish the same relation as in (\ref{eq:zot_matrix_trans2}). For
general first order terms of the form $\int_{\Omega} \psi_{i} {\bf b}
\cdot \nabla \phi_{j}$, with ${\bf b} \in [L^{\infty}(\Omega)]^{d}$,
it is simple to check that for the case the smaller element defines
the test space, the element matrix $M_{T'}$ can be calculated on the
smaller element and multiplied with ${\cal C}$ from the left. If the
smaller element defines the space of trial function, the element
matrix calculated on the smaller element must be multiplied with ${\cal
  C}_{\nabla}^{T}$ from the right.


\subsection{Calculation of transformation matrices}
We have shown that if the transformation matrix is calculated for a
given tuple of small and large element, the additional cost for
virtual mesh assembling is very small. In this section, we show how to
compute these transformation matrices efficiently. We formally define
a virtual element pair by the tuple
\begin{equation}
  \label{eq:vep}
  (T, \{\alpha_{0}, \text{\dots}, \alpha_{n}\}) = (T, \alpha) \: \text{with} \: \alpha_{i} \in \{\text{L}, \text{R}\},
\end{equation}
where $T$ is the larger element of the pair and $\alpha$ is an ordered
set that is interpreted as the refinement sequence for element
$T$. Thus, L denotes the ``left'', and R denotes the ``right''
children of and element. Furthermore, we define a function $TRA$ that
uniquely maps a virtual element pair to the smaller element. It is
defined recursively by:
\begin{eqnarray*}
  TRA(T, \emptyset) & = & T \\
  TRA(T, \{\alpha_{0}, \alpha_{1}, \text{\dots}, \alpha_{n}\}) & = &
  \begin{cases}
    TRA(T_{\text{L}}, \{\alpha_{1}, \text{\dots}, \alpha_{n}\}) & \text{if} \: \alpha_{0} = \text{L} \\
    TRA(T_{\text{R}}, \{\alpha_{1}, \text{\dots}, \alpha_{n}\}) & \text{if} \: \alpha_{0} = \text{R},
  \end{cases}
\end{eqnarray*}
where $T_{\text{L}}$ is the left child of the element $T$, and
$T_{\text{R}}$ the right child of the element, respectively. In the
same way we can now define transformation matrices as functions on
refinement sequences:
\begin{eqnarray*}
  {\cal C}(\emptyset) & = & \text{I} \\
  {\cal C}(\{\alpha_{0}, \alpha_{1}, \text{\dots}, \alpha_{n}\}) & = &
  \begin{cases}
    {\cal C}_{\text{L}} \cdot {\cal C}(\{\alpha_{1}, \text{\dots}, \alpha_{n}\}) & \text{if} \: \alpha_{0} = \text{L} \\
    {\cal C}_{\text{R}} \cdot {\cal C}(\{\alpha_{1}, \text{\dots}, \alpha_{n}\}) & \text{if} \: \alpha_{0} = \text{R},
  \end{cases}
\end{eqnarray*}
where ${\cal C}_{L}$ and ${\cal C}_{R}$ are the transformation matrix
for the left child and the right child, respectively, of the reference
element.

\subsection{Implementation issues}
Although the calculation of transformation matrices is quite fast, it
can considerably increase the time for assembling the linear
system. This is especially the case, if one mesh is much coarsen in
some regions than the other mesh. To circumvent this problem, we have
implemented a cache, that stores the transformation matrices. In the
mesh traverse routine, an 64 bit integer data type stores the
refinement sequence bit-wise, as it is defined by (\ref{eq:vep}). If
the bit on the i-th position is set, the finer element is a
right-refinement of its parent element, otherwise it is a
left-refinement of it. Of coarse this limits the level gap between the
coarser and the finer level to be less or equal to 64. But we have not
found any practical simulations, where this value is more than
30. Using this data type makes it than easy to define associative array
that uniquely maps a refinement sequence to a transformation
matrices. If a transformation matrix was computed for a given
refinement sequence for the first time, it will be stored in this
array. To access previously computed matrices using the integer key is
than very cheap. In general this data structure should be restricted to
a fixed number of matrices to not spend to much of memory. In all of
our simulations, the number of matrices that should be stored in the
cache never exceed 100.000. Also in the 3D case with linear elements
the overall memory usage is than around 2 Mbyte, and can thus be
neglected. Therefore, we have not yet considered to implement an upper
limit for the cache.

\section{Numerical results}
\label{sec:numerical_results}

In this section, we present several examples, where the multi-mesh
approach is superior in contrast to the standard single-mesh finite
element method. Examples to be considered are phase-field equations to
study solid-liquid and solid-solid phase transitions. For a recent
review we refer to \cite{ElderProvatas_Wiley_2010}. These equations
involve at least one variable, the phase-field variable, which is
almost constant in most parts of the domain and thus can be
discretized within these parts using a coarse mesh. Within the
interface region however a high resolution is required to resolve the
smooth transition between the different phases. A second variable in
such systems is typically a diffusion field which in most cases varies
smoothly across the whole domain and thus will require a finer mesh
outside of the interface and a coarser mesh within the interface. Such
problems are therefor well suited for our multi-mesh approach. We will
consider dendritic growth in solidification and coarsening phenomena
in binary alloys to demonstrate the applicability.

Other examples for which large computational savings due to the use of
the multi-mesh approach are expected are diffuse interface and diffuse
domain approximations for PDE's to be solved on surfaces are within
complicated domains. The approaches introduced in
\cite{RaetzVoigt_CMS_2006} and \cite{LiLowengrubRaetzVoigt_CMS_2009},
respectively, use a phase-field function to implicitly describe the
domain the PDE has to be solved on.  For the same reason as in
phase-transition problems the distinct solution behaviour of the
different variables will lead to large savings if the multi-mesh
approach is applied. The approach has already been used for
applications such as chaotic mixing in microfluidics
\cite{AlandLowengrubVoigt_CMES_2010}, tip splitting of droplets with
soluble surfactants \cite{TeigenLowengrubVoigt_JCP_2010}, and
chemotaxis of stem cells in 3D scaffolds
\cite{LandsbergStengerGelinskyRoesenWolffVoigt_JMB_2010}.

As a further example we demonstrate that the multi-mesh approach can
also be used to easily fulfill the inf-sup condition for saddle-point
problems if both components are discretized using linear Lagrange
elements. We demonstrate this numerically for the incompressible
Navier-Stokes equation with piecewise linear elements for velocity and
pressure, but a finer mesh used for the velocity. Such an approach
might be superior to mixed finite elements of higher order or
stabilized schemes in terms of computational efficiency and
implementational efforts. For a review on efficient finite element
methods for the incompressible Navier-Stokes equation we refer to
\cite{Turek_Springer_1999}. We demonstrate the applicability of the
multi-mesh approach on the classical driven cavity problem.

\subsection{Dendritic growth}
We first consider dendritic growth using a phase-field model, which
today is the method of choice to simulate microstructure evolution
during solidification. For reviews we refer to
\cite{BoettingerWarrenBeckermannKarma_RMR_2002}. A
widely used model for quantitative simulations of dendritic structures
was introduced by Karma and Rappel \cite{KarmaRappel_PRE_1996,
  Karma_1998}, which reads in non-dimensional from
\begin{equation}
  \begin{split}
    A^{2}(n) \partial_{t} \phi &= (\phi - \lambda u (1 - \phi^{2}))(1
    - \phi^{2}) + \nabla (A^{2}(n) \nabla \phi) +
    \sum_{i = 1}^{d} \partial_{x_{i}} \left( |\nabla \phi|^{2}A(n)\frac{\partial A(n)}{\partial_{x_{i}} \phi}\right)\\
    \partial_{t} u &=
    D \nabla^{2} u + \frac{1}{2} \partial_{t}\phi,
  \end{split}
\end{equation}
where $d = 2,3$ is the dimension, $D$ is the thermal diffusivity
constant, $\lambda = \frac{D}{a_{2}}$, with $a_{2} = 0.6267$ is a
coupling term between the phase-field variable $\phi$ and the thermal
field $u$ and $A$ is an anisotropy function. For both, simulation in
2D and 3D, we use the following anisotropy function:
\begin{equation}
   A(n) = 
  (1 - 3 \epsilon) 
  \left(1 + \frac{4 \epsilon}{1 - 3\epsilon} 
  \frac{\sum_{i = 1}^{d}\phi_{x_{i}}^{4}}{|\nabla \phi|^{4}}\right),
\end{equation}
where $\epsilon$ controls the strength of the anisotropy and $n =
\frac{\nabla \phi}{|\nabla \phi|}$ denotes the normal to the
solid-liquid interface. In this setting the phase-field variable is
$-1$ in liquid and $1$ in solid, and the melting temperature is set to
be zero. As boundary condition we set $u = - \Delta$ to specify an
undercooling. For the phase-field variable we use zero-flux boundary
conditions.

To implement these equations in our toolbox AMDiS, we first discretize
in time. This is here done using a semi-implicit Euler method, which
yields a sequence of nonlinear elliptic PDEs:
\begin{equation}
  \begin{split}
    \frac{A^{2}(n_n)}{\tau}\phi_{n+1} + f + g - \nabla (A^2(n_n)
    \nabla \phi_{n+1}) - {\cal L}[A(n_n)] &=
    \frac{A^{2}(n_n)}{\tau}\phi_{n} \\
    \frac{u_{n+1}}{\tau} - D \nabla^{2}u_{n+1} - \frac{1}{2}
    \frac{\phi_{n+1}}{\tau} &= \frac{u_{n}}{\tau} - \frac{1}{2}
    \frac{\phi_{n}}{\tau}.
  \end{split}
\end{equation}
with $f = \phi_{n+1}^3 - \phi_{n+1}$, $g = \lambda(1 - \phi_{n+1}^2)^2
u_{n+1}$ and 
\begin{equation*}
  {\cal L}[A(n_n)] = \sum_{i = 1}^{d} \partial_{x_{i}}
  \left( |\nabla \phi_{n+1}|^{2}A(n_{n})\frac{\partial A(n_{n})}{\partial_{x_{i}}
      \phi_{n}}\right).
\end{equation*}
We now linearize the involved nonlinear terms $f$ and $g$:
\begin{equation}
  \begin{split}
    f &\approx (3 \phi_{n}^{2} - 1) \phi_{n+1} - 2 \phi_{n}^{3} \\
    g &\approx \lambda(1 - \phi_{n}^{2})^{2}u_{n+1}
  \end{split}
\end{equation}
and obtain a linear system for $\phi_{n+1}$ and $u_{n+1}$ to be solved
in each time step.

To compare our multi-mesh method with a standard adaptive finite
element approach, we have computed a dendrite using linear finite
elements. The following parameters are used:
\begin{equation*}
  \Delta = 0.65, D = 1.0, \epsilon = 0.05
\end{equation*}
We have run the simulation with a constant timestep $\tau = 1.0$ up to
time $7500$. To speedup the computation we have employed the symmetry
of the solution and limited the computation to the upper right
quadrant with a domain size of 800 into each direction.  The adaptive
mesh refinement relies on the residuum based a posteriori error
estimate, as defined by \ref{eq:residual_ee}. By setting $C_{0}$ to 0
and $C_{1}$ to 1, we restrict the estimator to the jump residuum
only. We have set the tolerance to be $tol_{\phi} = 0.5$ and $tol_{u}
= 0.25$. For adaptivity, the equidistribution strategy with parameters
$\theta_{R} = 0.8$ and $\theta_{C} = 0.2$ was used. Thus, the
interface is resolved by around 20 points.

\begin{figure}[bth]
  \begin{center}
   \includegraphics[width=10cm]{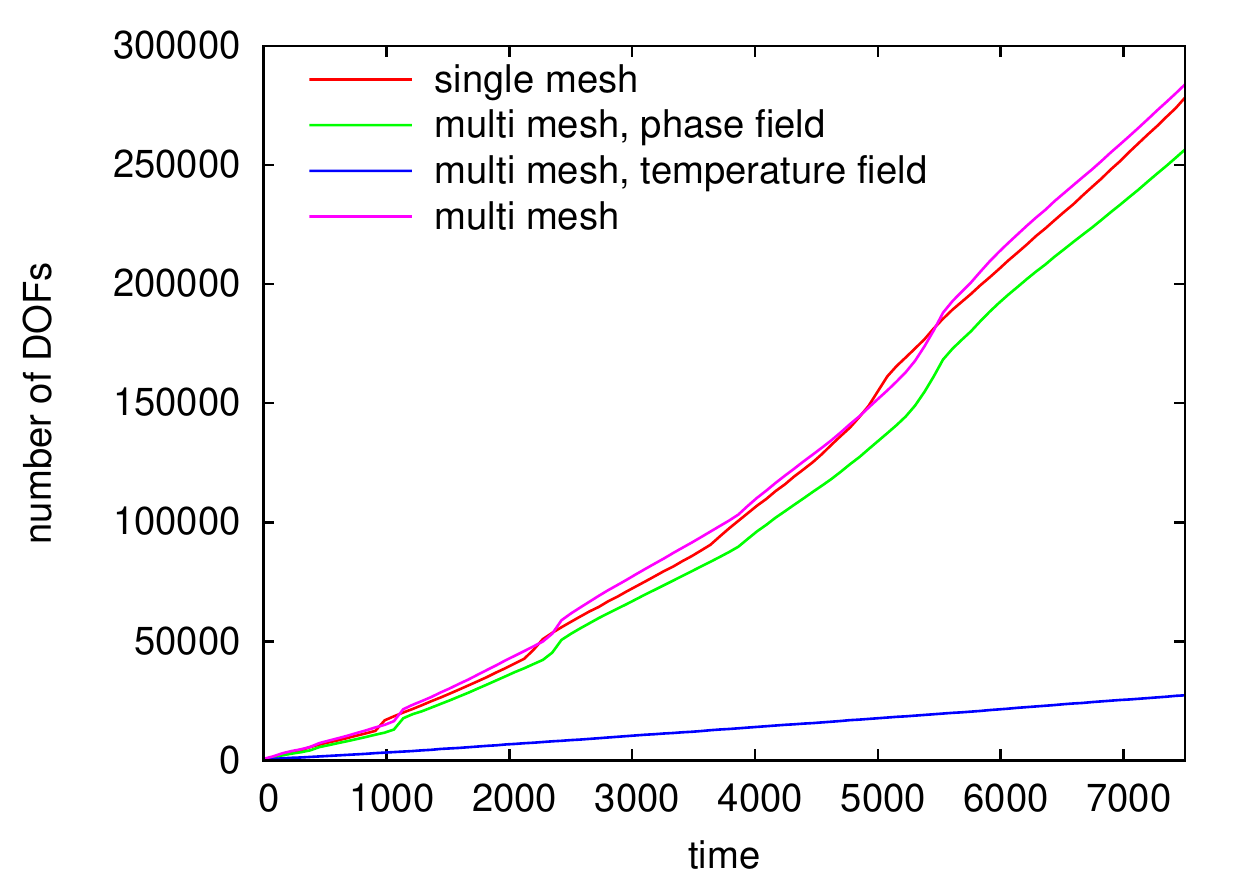}
   \caption{Evolution of degrees of freedom in time for the phase
     field and the temperature field using single-mesh and multi-mesh
     method.}
   \label{fig:ev_dof}
  \end{center}
\end{figure}

\begin{figure}[bth]
  \begin{center}    
    \includegraphics[width=\textwidth]{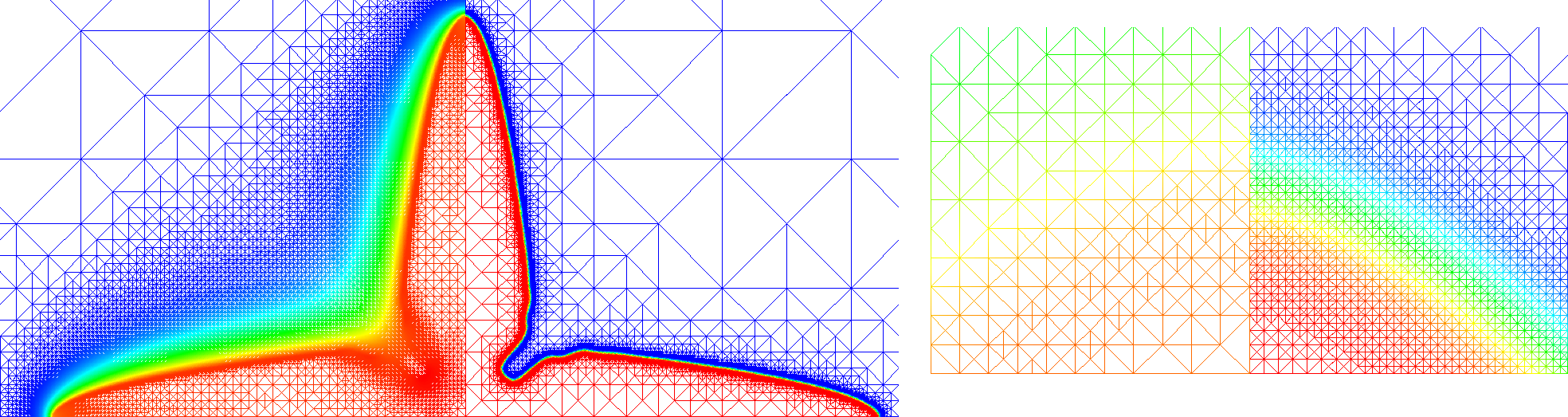}
    \caption{a) 2D dendrite computed at $t = 7500$ using the
      multi-mesh method with the parameters $\Delta = 0.65$, $D =
      1.0$, $\epsilon = 0.05$ and a timestep $\tau = 1.0$. Left shows
      the mesh of the temperature field and right shows the mesh of
      the phase field. b) Zoom into the upper tip showing the
      different resolution of both meshes.}
   \label{fig:dend2d_mesh}
  \end{center}
\end{figure}

The result of both computations coincides at the final timestep. As a
quantitative comparison we use the tip velocity of the dendrite. As
reported by Karma and Rappel \cite{Karma_1998}, for this parameter set
analytical calculations lead to a steady-state tip velocity $V_{tip} =
0.0469$. In both of our calculations, the tip velocity varies around
1\% to this value. Using the single-mesh method, within the final
timestep the mesh consist of 278.252 degrees of freedom (DOFs). When
our multi-mesh method is used, the same solution can be obtained with
a mesh for the phase-field with 256.342 DOFs and 27.354 DOFs for the
temperature field. The gap between these numbers increase over time,
see Figure \ref{fig:ev_dof} showing the development of DOFs over
time. Figure \ref{fig:dend2d_mesh} qualitatively compares the meshes
of the phase-field variable and the temperature field which shows the
expected finer resolution of the phase-field mesh within the interface
and its coarser resolution within the solid and liquid region.

The computational time for both methods is compared in Table
\ref{tab:compare_methods}. The assemblage procedure in the multi-mesh
methods is around 6\% faster, although computing the element matrices
is slower due to the extra matrix-matrix multiplication. This is
easily explained by the fact that we have much less element matrices
to compute and the overall matrix data structure is around 50\%
smaller with respect to the number of non-zero entries. This is also
reflected in the solution time for the linear system. We have run all
computation twice, with using either UMFPACK, a multifrontal sparse
direct solver \cite{Davis_2004}, or the BiCGStab($\ell$) with diagonal
preconditioning, that is part of the Math Template Library (MTL4), see
\cite{Gottschling_2007}. When using the first one within the
multi-mesh method, the solution time can be reduced by 40\% and also
the memory usage, which is often the most critical limitation in the
usage of direct solvers, is reduced in this magnitude. An even more
drastic reduction of the computation time can be achieved when using
an iterative solver. Here, the number of iterations is around 30\%
less with the multi-mesh method and each iteration is faster due to
the smaller matrix. The time for error estimation is halved, as
expected, since it scales linearly with the number of elements in the
mesh. Altogether, the time reduction is significant in all parts of
the finite element method for this example. In addition the approach
also leads to a drastic reduction in the memory usage.

\begin{table}[tbh]
  \begin{center}
    \begin{tabular}{|l|l|l||l|}
      \hline
      & \textbf{single-mesh} & \textbf{multi-mesh} & \textbf{speedup} \\\hline
      assembler & 5.98s & 5.62s & 6.0\% \\\hline
      solver: UMFPACK & 6.69s & 4.12s & 38.4\%\\\hline
      solver: BiCGStab$(\ell)$ & 14.26s & 4.28s & 69.9\% \\\hline
      estimator & 3.40s & 1.71s & 49.7\% \\\hline \hline
      overall with UMFPACK & 16.07s & 11.45s & 28.7\% \\\hline
      overall with BiCGStab$(\ell)$ & 23.64s & 11.61s & 50.9\% \\\hline
    \end{tabular}
    \caption{Comparison of runtime when using single-mesh and
      multi-mesh method. The values are the average of 7500
      timesteps.}
    \label{tab:compare_methods}   
  \end{center}
\end{table}

\begin{figure}[bth]
  \begin{center}
   \includegraphics[width=\textwidth]{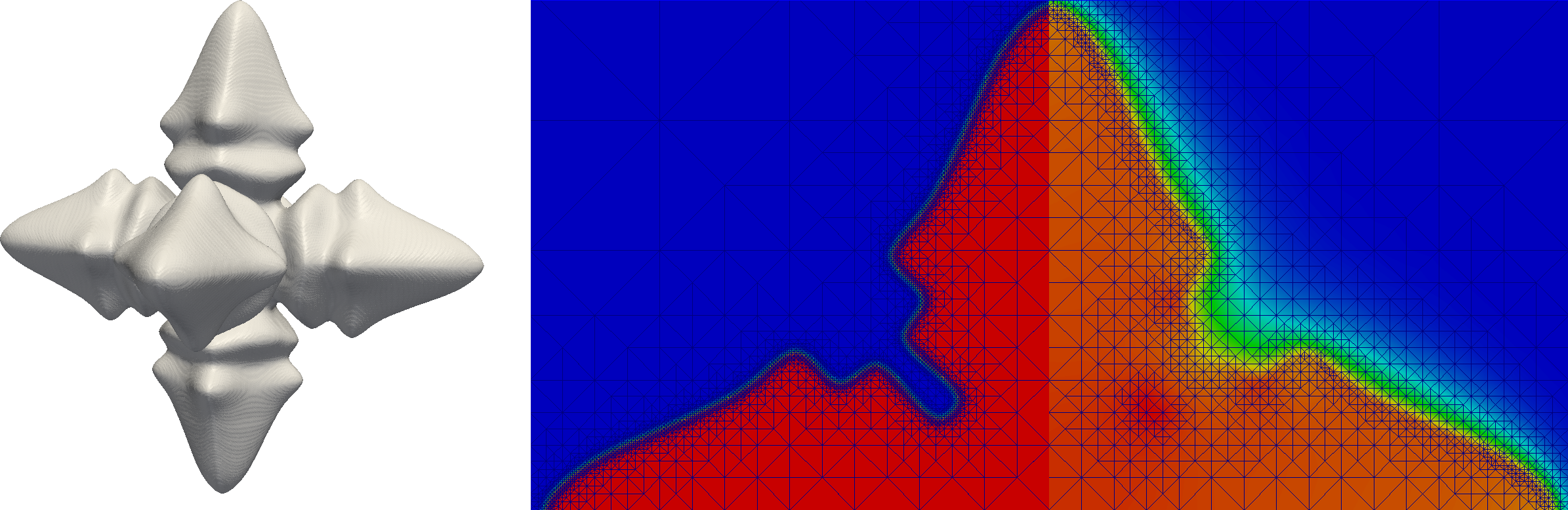}
   \caption{a) 3D dendrites at t = 2500 using the multi-mesh method
     with the parameters $\Delta = 0.55$, $D = 1$, $\epsilon = 0.05$
     and a timestep $\tau = 1.0$. b) Slice through the dendrites
     showing the mesh of the phase field on the left and the mesh of
     the temperature field on the right hand side.}
   \label{fig:dendrites3d}
  \end{center}
\end{figure}

The results are even more significant in 3D. Figure
\ref{fig:dendrites3d} shows the result of computing a dendrite with
the multi-mesh method and the following parameters:
\begin{equation*}
  \Delta = 0.55, D = 1, \epsilon = 0.05
\end{equation*}
We have run the simulation with a constant timestep $\tau = 1.0$ up to
time 2500. The evolution of degrees of freedom over time is quite
similar to the 2D example. When using the multi-mesh method, the time
for solving the linear system, again using the BiCGStab($\ell$) solver
with diagonal preconditioning, can be reduced by around 60\%. The time
for error estimation is around half the time needed by the single-mesh
method. Because the time for assembling the linear system is more
significant in 3D, the overall time reduction with the multi-mesh
method is around $24.4\%$.

\subsection{Coarsening}

As a second example we consider coarsening of a binary structure using
a Cahn-Hilliard equation. We here concentrate only on the
phenomenological behaviour of the solution and thus consider the
simplest isotropic model, which reads
\begin{equation}
  \begin{split}
    \partial_t \phi &= \Delta \mu \\
    \mu &= - \epsilon \Delta \phi + \frac{1}{\epsilon} G^\prime(\phi)
  \end{split}
\end{equation}
for a phase-field function $\phi$ and a chemical potential $\mu$. The
parameter $\epsilon$ again defines a length scale over which the
interface is smeared out, and $G(\phi) = 18 \phi^2(1- \phi)^2$ defines
a double-well potential. To discretize in time we again use a
semi-implicit Euler scheme
\begin{equation}
  \begin{split}
    \frac{1}{\tau} \phi_{n+1} - \Delta \mu_{n+1} &= \frac{1}{\tau} \phi_n \\
    \mu_{n+1} + \epsilon \Delta \phi_{n+1} - \frac{1}{\epsilon} G^\prime(\phi_{n+1}) &= 0
  \end{split}
\end{equation}
in which we linearize $G^\prime(\phi_{n+1}) \approx
G^{\prime\prime}(\phi_{n})\phi_{n+1} + G^\prime(\phi_n) -
G^{\prime\prime}(\phi_n) \phi_n$.

To compare our multi-mesh method with a standard adaptive finite
element approach, we have computed the spinodal decomposition and
coarsening process using Lagrange finite elements of fourth order. We
use $\epsilon = 5\cdot 10^{-4}$. The adaptive mesh refinement relies
on the residuum based a posteriori error estimate. As we have done it
in the dendritic growth simulation, also here only the jump residuum
is considered, i.e., the constants $C_{0}$ and $C_{1}$ in Equation
(\ref{eq:residual_ee}) are set to 0 and 1, respectively. For both
methods, the error tolerance are set to $tol_{\phi} = 2.5 \cdot
10^{-4}$ and $tol_{\mu} = 5 \cdot 10^{-2}$. For adaptivity, the
equidistribution strategy with parameters $\theta_{R} = 0.8$ and
$\theta_{C} = 0.2$ was used. Using these parameters, the interface is
resolved by around 10 points.

\begin{figure}[bth]
  \begin{center}
    \includegraphics[width=\textwidth]{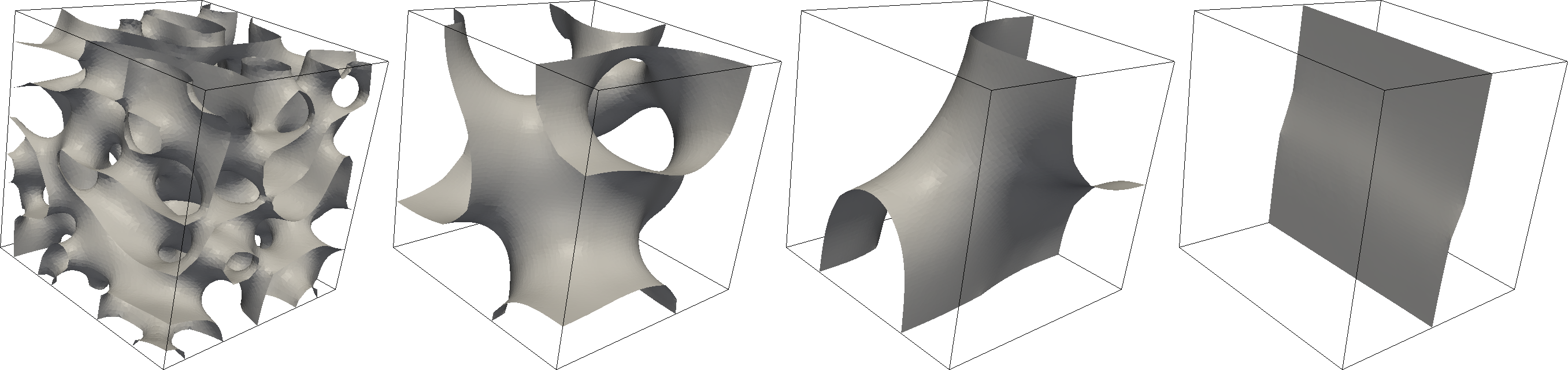}
    \caption{Solution of the Cahn-Hilliard equation for $t = 0.02$, $t
      = 1.0$, $t = 5.0$ and $t = 12.50$}
    \label{fig:cahn_hilliard_3d}
  \end{center}
\end{figure}

\begin{figure}[bth]
  \begin{center}
    \includegraphics[width=10cm]{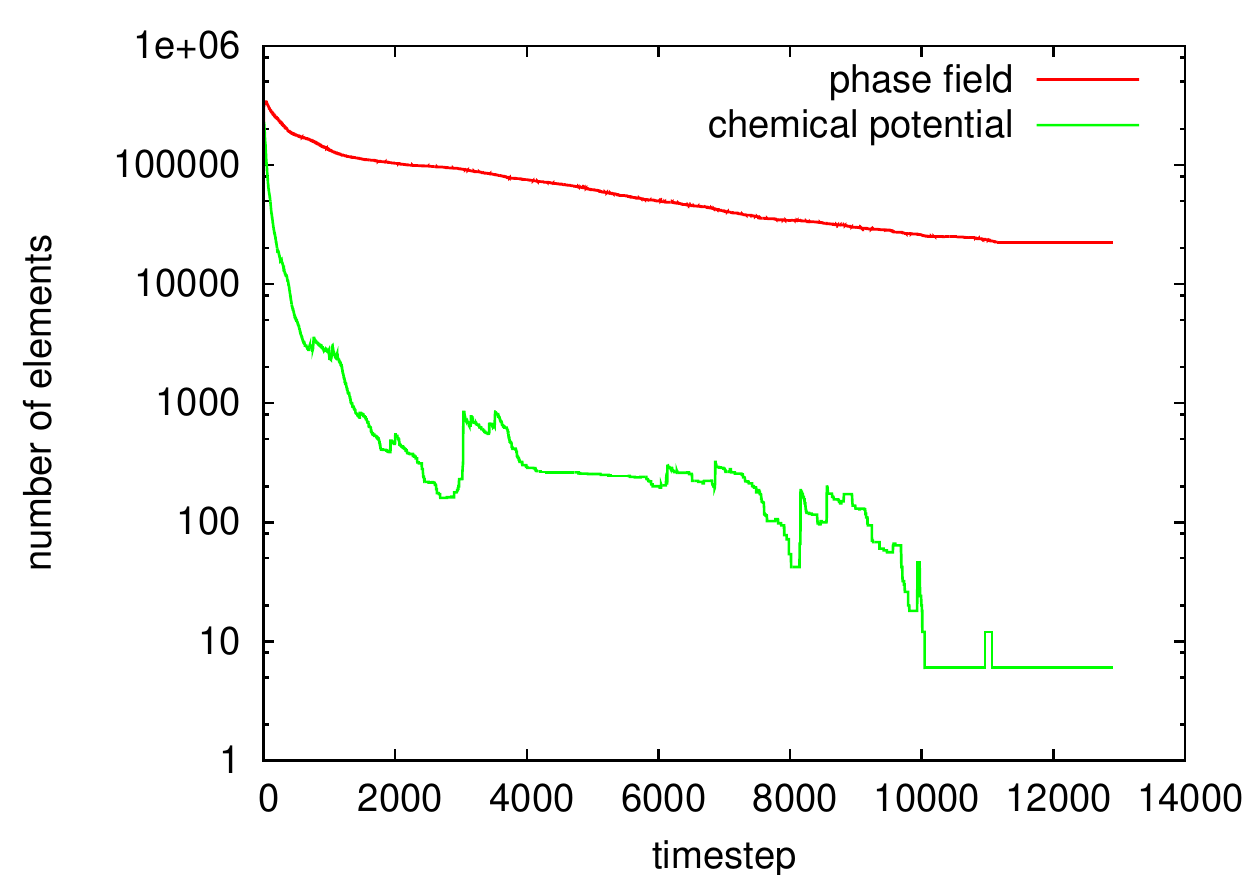}
    \caption{Evolution of number of elements for both variables of the
      Cahn-Hilliard equation.}
  \end{center}
  \label{fig:cahn_hilliard_elems}
\end{figure}

The simulation was started from noise. The first mesh was globally
refined with 196608 elements. The constant timestep was chosen to be
$\tau = 10^{-3}$. We have disabled the adaptivity for the first 10
timesteps, until a very first coarsening in the domain was
achieved. Then the simulation was executed up to $t = 13.0$, where
both phases are nearly separated. Figure \ref{fig:cahn_hilliard_3d}
shows the phase field, i.e., the 0.5 contour of the first solution
variable, for four different timesteps. The number of elements and
degrees of freedom is linear to the area of the interface that must be
resolved on the domain. Indeed, the chemical potential can be resolved
on a mush coarser grid, since it is independent of the resolution of
the phase field. In the final state, the chemical potential is
constant on the whole domain, and the macro mesh (which consists of 6
elements in this simulation) is enough to resolve it. The evolution of
the number of elements for both variable over time is plotted in
Figure \ref{fig:cahn_hilliard_elems}. As expected, the number of
elements for the phase field monotonously decreases as its area
shrinks due to the coarsening process. The number of elements for the
chemical potential rapidly decreases at the very first beginning, as
the initial mesh is over refined to resolve this variable. For most of
the simulation, the number of elements of the chemical potential is
three magnitudes smaller the the number of elements for the phase
field variable. This gap is also reflected in the computation time for
the single-mesh and the multi-mesh method. AMDiS has computed the
solution with the multi-mesh method in 722 minutes. When using the
standard single-mesh method, where the number of elements for the only
mesh is nearly equivalent to the number of elements needed in the
phase field mesh when using the multi-mesh approach, the simulation
time was 1465 minutes. For both simulations, we have used the
BiCGStab$(\ell)$ solver with diagonal preconditioning.

\subsection{Fluid dynamics}
For the last example, we use our multi-mesh method to solve problems
in fluid dynamics using standard linear finite elements. The inf-sub
stability is established by using two different meshes. In 2D, the
mesh for the velocity components is refined twice more than the mesh
for pressure. In the 3D case, the velocity mesh has to be refined
three times to get the corresponding refinement structure. This
discretization was introduced by Bercovier and Pironneau
\cite{Bercovier_1979}, and was analyzed and proven to be stable by
Verf\"urth \cite{Verfurth_1984}. To show this technique, we solve the
standard instationary Navier-Stokes equation given by
\begin{equation}
  \label{eq:ns}
  \begin{split}
    \partial_{t} u -\nu \nabla^{2} u + u \cdot \nabla u + \nabla p &= f \\ 
    \nabla \cdot u &= 0,
  \end{split}
\end{equation}
where $\nu > 0$ is the kinematic viscosity. The time is discretized by
the standard backward Euler method. The nonlinear term in
(\ref{eq:ns}) is linearized by
\begin{equation}
  u_{n} \cdot \nabla u_{n+1}
\end{equation}

\begin{figure}[bth]
   \begin{center}
    \begin{tabular}{lcc}
      single mesh &
      \includegraphics[width=4.0cm]{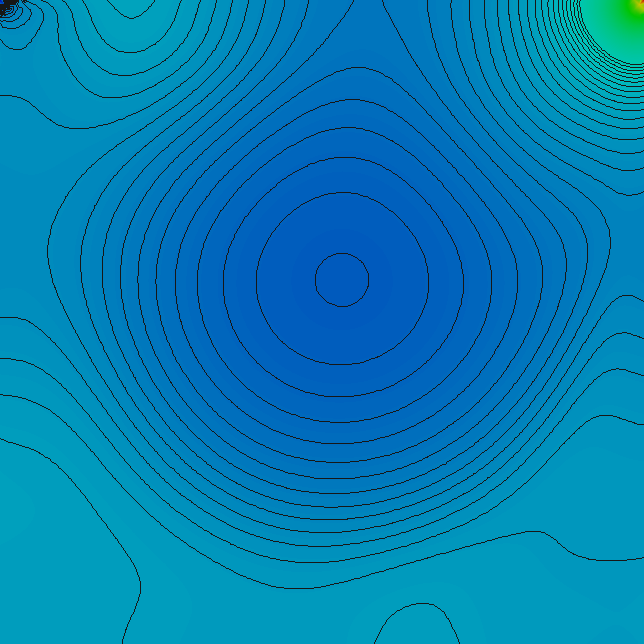} &
      \includegraphics[width=4.0cm]{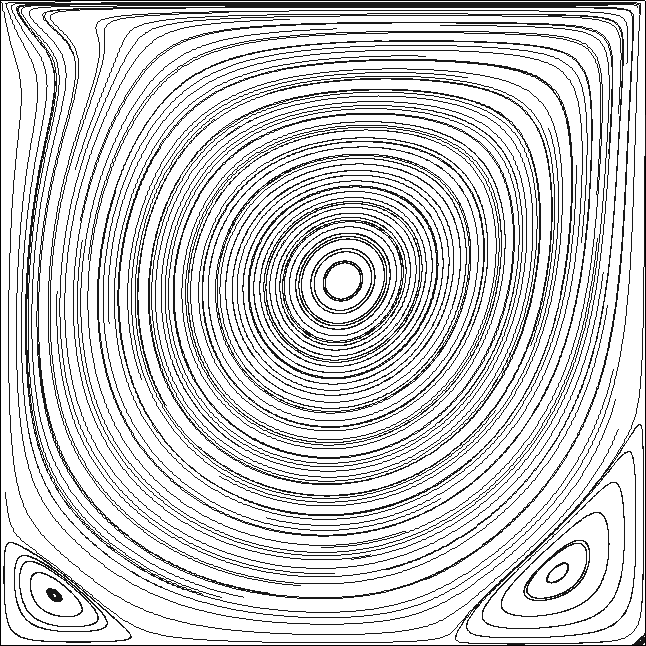} \\
      multi mesh &
      \includegraphics[width=4.0cm]{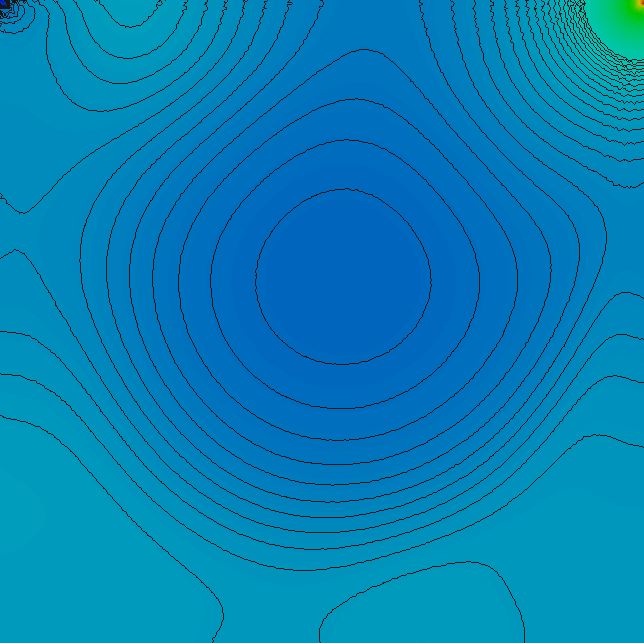} &
      \includegraphics[width=4.0cm]{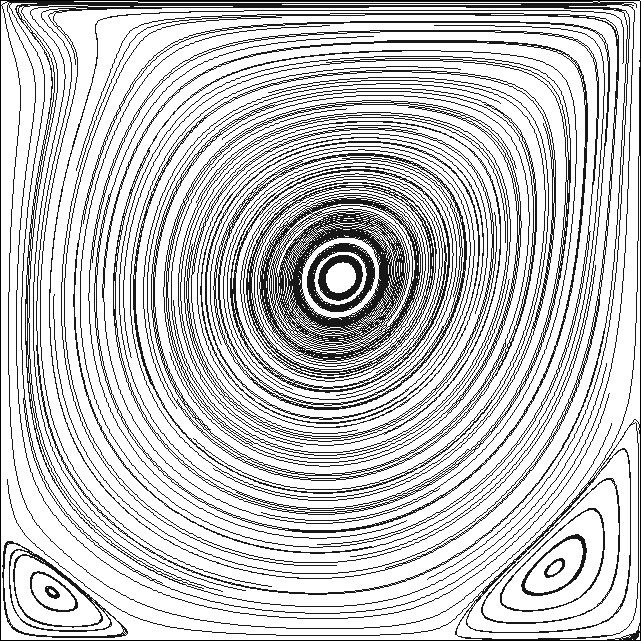}
    \end{tabular}
  \end{center}

  \label{fig:driven_cavity_re1000}
  \caption{Results for $Re = 1000$. The left column shows 30 pressure
    isolines in the range $-0.12,0.12$. The right column shows the
    streamlines contours of the velocity.}
\end{figure}

The model problem is the ``driven cavity'' flow, as described and
analyzed in \cite{Wall_1999, Ghia_1982}. In a unit square, the
boundary conditions for the velocity are set to be zero on the left,
right and lower part of the domain. On the top, the velocity into x
direction is set to be one and into y direction to be zero. In the
upper corners, both velocities are set to be zero, which models the so
called non-leaky boundary conditions. The computation was done for
several Reynold numbers varying between 50 and 1000. First, we have
used the single-mesh method with a standard Taylor-Hood element, i.e.,
second order Lagrange finite element for the velocity components a
linear Lagrange finite element for the pressure. Than, we have
compared these results with our multi-mesh method, where for both
components a linear finite elements was used and, instead, the mesh
for velocity discretization is refine twice more than the
pressure. All computations were done with a fixed timestep $\tau =
0.01$ and aborted, when the relative change in velocity and pressure
is less than $10^{-6}$. Figure \ref{fig:driven_cavity_re1000} shows
the results for $Re = 1000$. In Table \ref{tab:compare_driven_cavity},
we give the position of all eddies and compare our results with the
work of Ghia et al.\ \cite{Ghia_1982} and Wall \cite{Wall_1999}.

\begin{table}[tbh]
  \begin{center}
    \begin{tabular}{|l|l|l|l|l|}
      \hline
      & Eddy 1 & Eddy 2 & Eddy 3 & Eddy 4 \\\hline
      single-mesh & 0.5310, & 0.8633, & 0.0838, & 0.9937, \\
      & 0.5658 & 0.1116 & 0.0775 & 0.0062 \\\hline
      multi-mesh & 0.5305, & 0.8669, & 0.0813, & 0.9953, \\
      & 0.5671 & 0.1125 & 0.0750 & 0.0062 \\\hline
      Wall & 0.5308, & 0.8643, & 0.0832, & 0.9941, \\
      & 0.5660 & 0.1115 & 0.0775 & 0.0066 \\\hline
      Ghia et al.\ & 0.5313, & 0.8594, & 0.0859, & 0.9922, \\
      & 0.5625 & 0.1094 & 0.0781 & 0.0078 \\\hline
    \end{tabular}
    \caption{Comparison of eddy position for $Re = 1000$ in the driven
      cavity model.}
    \label{tab:compare_driven_cavity}   
  \end{center}
\end{table}

In both, the single-mesh method and the multi-mesh method, all finite
element spaces have the same number of unknowns. This is the reason,
why the usage of the latter one is not faster in contrast to the
single-mesh method, as it is the case in dendritic growth. The time
for assembling the linear system growth from $4.13$ seconds to $5.79$
seconds, which is mainly caused by the multiplication of the element
matrices with the transformation matrices. Instead, the average
solution time with a BiCGStab$(\ell)$ solver and ILU preconditioning
decreases from $10.18$ seconds to $8.88$ seconds. Although the linear
systems have the same number of unknowns, the linear systems resulting
from the single-mesh method are denser due to the usage of second
order finite elements. The number of non-zero entries decreases around
20\% when linear elements are used on both meshes.

\section{Conclusion}

To further improve efficiency of adaptive finite element simulations
we consider the usage of different adaptively refined meshes for
different variables in systems of nonlinear, time-depended PDEs. The
different variable can have very distinct solution behaviour. To
resolve this the meshes can be independently adapted for each
variable. Our multi-mesh method works for Lagrange finite elements of
arbitrary degree and is independent of the spatial dimension.  The
approach is well defined, and can be implemented in existing adaptive
finite element codes with minimal effort. We have demonstrated for
various examples that the resulting linear systems are usually much
smaller, when compared to the usage of a single mesh, and the overall
computational runtime can be more than halved in various cases.  Phase
transition problems within a diffuse interface approach are well
suited for our approach. The same holds for saddle-point problems in
which the inf-sup condition can be fullfilled for finite elements of
the same order.

Further examples which are currently under investigation include
general diffuse interface concepts to solve PDEs in complex
domain. Here a phase-field function is used to describe the domain
implicitly \cite{LiLowengrubRaetzVoigt_CMS_2009}, which only requires
a fine resolution along the boundaries. The approach might also be
used in time stepping schemes to prevent loss of information during
coarsening. In a classical approach the solution from the old time
step is simply interpolated to the new mesh at the new time step. If
the new mesh is coarser information is lost, which can be prevented by
using the multi-mesh approach for the solution at different time
steps. And also in optimal control problems the approach is very
promizing, as in many situations the dual solution is much smoother
than the primal solution and thus can be discretized on a much coarser
mesh using are multi-mesh approach which will be demonstrated for
control of an Allen-Cahn equation.

\section{Acknowledgement}

We would like to thank Rainer Backofen for fruitful discussions. The
work has been supported by DFG through Vo899/5-1 and Vo899/11-1.

\bibliographystyle{unsrt}
\bibliography{fem}

\end{document}